\documentclass[a4paper, 12pt]{article}
\usepackage{enumerate, theorem}
\usepackage{amsmath, amsfonts, amssymb}
\usepackage[height=22.5cm, width=15cm]{geometry}
\usepackage[all]{xy}
\usepackage{mathrsfs, yfonts}
\usepackage[utf8]{inputenc}

\DeclareMathOperator{\id}{id}

\newcommand{\parag}[1]{\paragraph{\sc{#1.}}}

\newtheorem{thm}{Theorem}[subsection]

\newtheorem{cor}[thm]{Corollary}
\newtheorem{prop}[thm]{Proposition}
\newtheorem{lemma}[thm]{Lemma}

\begin{document}
\title{Complement to Higher Bernstein Polynomials and Multiple Poles of $\frac{1}{\Gamma(\lambda)}\int_X \vert f\vert^{2\lambda}\bar f^{-h}\rho\omega\wedge \bar \omega'$}

\author{Daniel Barlet\footnote{Barlet Daniel, Institut Elie Cartan UMR 7502  \newline
Universit\'e de Lorraine, CNRS, INRIA  et  Institut Universitaire de France, \newline
BP 239 - F - 54506 Vandoeuvre-l\`es-Nancy Cedex.France. \newline
e-mail : daniel.barlet@univ-lorraine.fr}.}

\maketitle

\parag{Abstract} We give, using higher Bernstein polynomials defined in our paper \cite{[B.23]}, a stronger version of our previous result in \cite{[B.22]} whose converse is proved in \cite{[B.23]} and we give some complements to the results in  \cite{[B.23]} which help to compute these higher order Bernstein polynomials. Then we show some non trivial  examples where we determine the root of the second Bernstein polynomial which is not a  double root of the full Bernstein polynomial and where the main theorem of \cite{[B.23]} applies  and localizes where a double pole exists for the meromorphic extension
of the (conjugate) analytic functional given by polar parts of  $\omega' \mapsto \vert f\vert^{2\lambda} \bar f^{-h}\rho \omega\wedge \bar \omega'$ when $h \in \mathbb{N}$ is large enough.

\parag{Classification AMS}  32 S 25; 32 S 40 ; 34 E 05 

\tableofcontents

\section{Introduction}

The first motivation for this ``complement" to the paper \cite{[B.23]} is to give some explicit examples of higher Bernstein polynomials in the case where they have roots which are not multiple roots of the (full) Bernstein polynomial of the fresco under consideration in {\it loc. cit.}, for instance in Theorem 8.5.1.\\
The first remark is that, in general, it is rather difficult to compute the Bernstein polynomial of the fresco associated to a given pair $(f, \omega)$, even in the case where $f$ has an isolated singularity. Nevertheless, in the case where $f$ is a polynomial in $\mathbb{C}[x_0, \dots, x_n]$ having $(n+2)$ monomials, we describe  in the article  \cite{[B.22]} a rather elementary method to obtain an estimation for the Bernstein polynomial of the fresco $\mathcal{F}_{f,\, \omega}$ associated to a monomial $(n+1)$-form $\omega$.\\
Of course, when the full Bernstein polynomial has a  root of multiplicity $k \geq 2$ then this root is also a root of the $j$-th Bernstein polynomial for each $j \in [1,k]$ but when the Bernstein polynomial has only simple roots, the computation of the higher order Bernstein polynomials, even in the special situation of   \cite{[B.22]}, is not obvious. We present here some examples where we show that the second Bernstein polynomial  is not trivial in cases where the full  Bernstein polynomial has no multiple root.\\
We concentrate here on the case of the function $f(x,y,z) = xy^3 + yz^3 + zx^3 + \lambda xyz $ where $\lambda \not= 0$ is any complex parameter, and where $\omega = \mu dx\wedge dy\wedge dz$ with 
$\mu$ a monomial in $\mathbb{C}[x,y,z]$. The tools used to estimate the Bernstein polynomial of the associated fresco $\mathcal{F}_{f, \omega}$ are valid for almost all polynomials in $\mathbb{C}[x_0, \dots, x_n]$ having $(n+2)$ monomials, for any $n \geq 1$ (see the conditions $C1$ and $C2$ in \cite{[B.22]}).\\
The tools used to determine the second Bernstein polynomial in the examples below are also easy to generalize, at least to obtain information on the smallest root of the Bernstein polynomial in the previous setting, knowing that this root is, in general, a root of the $k$-th Bernstein polynomial of the fresco $\mathcal{F}_{f, \omega}$ where $k$ is the nilpotent order of this fresco.\\
It seems not difficult, for some colleague fun of computers, to make a program which produces in such a situation, not only estimates for the (full) Bernstein polynomial of $\mathcal{F}_{f, \omega}$ but also, computing also the polynomial in (a,b) annihilating $[\omega]$ in $\mathcal{F}_{f, \omega}$ (described in \cite{[B.22]}) to estimate the smallest root of the Bernstein polynomial.\\

During the computations which lead to the examples presented in Section 4 I realized that some rather easy consequence of the Section 5 in \cite{[B.23]} where missing to enlighten the relationship between asymptotic expansions, themes and higher Bernstein polynomials. These results are given in Section 2 and will be add to the second version of my paper \cite{[B.23]} (in preparation).\\

Another point appears also during this period; the fact that the main result of \cite{[B.23]} is a converse of a statement which is more precise than the main result in \cite{[B.22]}. The reason is that, without the notion of higher Bernstein polynomial of a fresco, this stronger statement cannot be formulate ! So we add this improved version of the  main theorem of \cite{[B.22]} in Section 3 below.\\

\section{Some useful results to complete \cite{[B.23]}}

We begin by recalling some facts from \cite{[B.23]}. For the definition of a  fresco see Definition 5.1.2 or Section 7.1 in  \cite{[B.23]}.

 \parag{Reminder}\begin{enumerate}
 \item Any $[\alpha]$-primitive fresco $\mathcal{F}$ has an embedding in some $\Xi^{(N)}_\alpha \otimes V$  thanks to Theorem 5.1.3 in \cite{[B.23]}.
 \item Remind that, if $\mathcal{F}$ is any fresco, for each of its  higher Bernstein polynomial, $ B^k_{\mathcal{F}}$ for $k \geq 1$, we have 
 $$ B^k_{\mathcal{F}} = \prod_{\alpha} B^k_{\mathcal{F}^\alpha}$$
 where $\mathcal{F}^\alpha := \mathcal{F}\big/\mathcal{F}_{[\not= [\alpha]}$ is the $[\alpha]$-primitive quotient of $\mathcal{F}$.
 \item If a fresco has the root $-\alpha-m$ for its $k$-th Bernstein polynomial, where $k := d(\mathcal{F}^\alpha)$, there exists a $[\alpha]$-primitive rank $k$ quotient theme $T$ of $\mathcal{F}$
 such that the $k$-th Bernstein polynomial of $T$ is $(x+\alpha+m)$.
 \item Conversely, if $T$ is a $[\alpha]$-primitive quotient theme of rank $k$ of a $[\alpha]$-primitive fresco  $\mathcal{F}$ such that the $k$-th Bernstein polynomial of $T$ is $(x+\alpha+m)$, then there exists a root in 
 $[-\alpha-m, -\alpha] \cap \{\alpha -\mathbb{N}\}$ for the $k$-th Bernstein polynomial of $\mathcal{F}$. Moreover, when $d(\mathcal{F}) = k$, $(-\alpha -m)$ is a root of the  $k$-th Bernstein polynomial of $\mathcal{F}$
 \end{enumerate}

We give now a lemma and a remark which will be added in the second version of \cite{[B.23]} (in preparation).\\

\begin{lemma}\label{missing}
Let $\mathcal{F}$ be a $[\alpha]$-primitive fresco and assume that $-\alpha -m$ is a root of its $k$-th Bernstein polynomial. Then for each $j \in [1,k]$ there exists an integer $m_j \in [0, m]$ such that $-\alpha - m_j$ is a root of the $j$-th Bernstein polynomial of $\mathcal{F}$.
\end{lemma}

\parag{Proof}  By definition, if the nilpotent order for $\mathcal{F}$ is strictly bigger than $k$ then  the Bernstein polynomial of $S_k(\mathcal{F})$ has a root which is strictly bigger than $-\alpha -m$. So it is enough to prove the lemme when $k$ is the nilpotent order of $\mathcal{F}$. Then, by a descendant induction on $j \in [1, k-1]$ it is enough to prove the case $j = k-1$.\\
Taking the quotient by $S_{k-2}(\mathcal{F})$ we reduce the question in the case $k = 2$.\\
In this case,  there exists a quotient theme $T$ with rank $2$ whose Bernstein polynomial has the root $-\alpha-m$ as its minimal root. Then the other root $-\alpha-m'$ of the Bernstein polynomial of $T$ satisfies $-\alpha -m' \geq -\alpha -m$. Since $T$ is a quotient of $\mathcal{F}$, $-\alpha-m'$ is a root of the Bernstein polynomial of $\mathcal{F}$. If it is a root of the first Bernstein polynomial of $\mathcal{F}$ we are done. If this is not the case, $-\alpha-m'$ is a root of the second Bernstein polynomial of $\mathcal{F}$. But in this case $m' < m$ since the roots of the second Bernstein polynomial of $\mathcal{F}$ are simple. Then we can play the same game as before, but with the root $-\alpha -m'$. Since there is only finitely integer in $[0,m]$ we finally reach a root $-\alpha-m"$ of the first Bernstein polynomial of $\mathcal{F}$ such that $m" $ is in $[0, m]$.$\hfill \blacksquare$\\

 \parag{Remark} Let $T$ be a $[\alpha]$-primitive theme with rank $k$. Then its $j$-th Bernstein polynomial has degree $1$ for $j \in [1,k]$ and is equal to $(x+\alpha + m_j)$ where $-\alpha -m_1, \dots, -\alpha -m_k$ are the roots of its Bernstein polynomial  in decreasing order.\\
  So $-\alpha - m_k$ is the smallest root of its Bernstein polynomial.\\
  
  Note that, for any $[\alpha]$-primitive fresco, the smallest root of the Bernstein polynomial is always a root of the $k$-th Bernstein polynomial where $k = d(\mathcal{F})$ is the nilpotent order of the fresco $\mathcal{F}$. But  for a ``general" $[\alpha]$-primitive fresco, we do not know other  relation between the order of the roots of the Bernstein polynomial of $\mathcal{F}$ and the  roots of the $j$-th Bernstein polynomial  of $\mathcal{F}$ than the fact, given by the Lemma above.

  \begin{prop}\label{new 1}
  Let $\mathcal{F}$ be a semi-simple fresco. Then $-\lambda$ is a root of the Bernstein polynomial of $\mathcal{F}$ if and only if there exists a $\tilde{\mathcal{A}}$-linear surjective map
  $$\pi : \mathcal{F} \to \mathcal{E}_\lambda \simeq \tilde{\mathcal{A}}\big/ \tilde{\mathcal{A}}(a -\lambda b).$$
  \end{prop}
  
  \parag{Proof} The existence of $\pi$ is sufficient because the Bernstein polynomial of a quotient of $\mathcal{F}$  divides the the Bernstein polynomial of $\mathcal{F}$.\\
  Conversely, if $\lambda$ is a root of the Bernstein polynomial of $\mathcal{F}$, since $\mathcal{F}$ is semi-simple, there exists a Jordan-H\"older sequence for $\mathcal{F}$ such its last quotient is $\mathcal{E}_\lambda$, thanks to Proposition 7.2.1  in \cite{[B.23]}. So the proof is complete.$\hfill \blacksquare$.
  
  \begin{cor}\label{new 2}
  Let  $\mathcal{F}$ be be a $[\alpha]$-primitive fresco with nilpotent order $k$. Assume that $\mathcal{F} = \tilde{\mathcal{A}}e \subset \Xi^{(k-1)}_\alpha \otimes V$\footnote{This is not restrictive thanks to the results of Section 5 in \cite{[B.23]}.}. Let $p$ be the rank of $\mathcal{F}\big/S_{k-1}(\mathcal{F})$. Then there exists $p$ linearly  independent vectors $v_1, \dots, v_p$ in $V$ such that $e$ may be written
  $$ e = \sum_{j=1}^p S_j(b)s^{\alpha + m_j -1}(Log\, s)^{k-1}\otimes v_j  + \psi. $$
  where $\psi$ is in $\Xi^{(k-2)}_\alpha \otimes V$, and  where the  $S_j$ are invertible elements in $\mathbb{C}[[b]]$. Moreover we may  choose the vectors $v_1, \dots, v_p$ such  that $m_1 < \dots < m_p$.\\
  When this condition is fulfilled  the $k$-th Bernstein polynomial of $\mathcal{F}$ is  equal to $\prod_{j=1}^p(x + \alpha + m_j)$.
  \end{cor}
  
  For $\alpha = 1$ it is convenient to replace $\Xi^{(k-1)}_1$ by $\Xi^{(k)}_1\big/ \Xi^{(0)}_1$ to consider only the singular part of the asymptotic expansions. This is the case in the examples computed in Section $4$.
  
  \parag{Proof} Since $S_{k-1}(\mathcal{F}) = \mathcal{F} \cap (\Xi^{(k-2)}_\alpha \otimes V)$, it is enough to treat the semi-simple case. In this case, since each $\mathcal{E}_{\alpha+m}$ is embedded in $\Xi^{(0)}_\alpha$ we may assume that $e = \sum_{j=1}^q S_j(b)s^{\alpha+\mu_j-1}\otimes v_j $ where $v_1, \dots, v_q$ is a basis of  $V$ (by definition of semi-simplicity), where $m_1, \dots, m_q$ are non negative integers and where $S_j$ are invertible elements in $\mathbb{C}[[b]]$ or vanish identically. Moreover, since the saturation $\mathcal{F}^\sharp$ is a direct sum of $\mathcal{E}_{\alpha + m}$ and has the same rank than $\mathcal{F}$, we may assume that  the vector $v_j$ for which $S_ j \not= 0$ generate a subspace $W$ of dimension $p$ in $V$, where $p$ is the rank of $\mathcal{F}$. \\
  If the integer $\mu_1, \dots, \mu_p$ are pairwise distinct we may order the $v_1, \dots, v_p$ such that $\mu_1 < \dots < \mu_p$ and put $m_j:= \mu_j$. If this is not the case, consider $m_1$ the infimum of the $\mu_j$ and when $\mu_j = m_1$  let $w_1 = v_1 + \sum c_jv_j$ where the sum is on  each  $j \geq 2$ such that $\mu_j = m_1$ and where $c_j = S_j(0)S_1(0)^{-1}$ with $\mu_j = m_1$. 
   Now we obtain a new expression for $e$  in the basis $w_1, v_2, \dots, v_p$ of $W$, where $m_1$ is strictly less than all $\mu'_j$ which appear for $j \geq 2$. Continuing in this way we obtain that $(w_1, \dots, w_p)$ is a new basis of $W$ and $m_1 < \dots < m_p$. \\
   Then consider the $\tilde{\mathcal{A}}$-linear maps given by the linear forms $l_j \in V^*$ defined by $l_j(w_h) = \delta_{j,h}, h \in [1, p]$. The $\tilde{\mathcal{A}}$-linear map $\id\otimes l_j$ for $j \in [1,p]$ sends surjectively $\mathcal{F}$ to $\mathcal{E}_{\alpha + m_j}$ and this implies that $-(\alpha + m_j)$ is a root of the Bernstein polynomial of $\mathcal{F}$ for each $j \in [1, p]$. But since $\mathcal{F}$ has rank $p$ we obtain all the roots of its Bernstein polynomial since the $m_j$ are pair-wise distinct\footnote{Note that the initial $\mu_j$ gives also roots of the $k$-th Bernstein polynomial of $\mathcal{F}$ but they may not give all the roots.}. This completes the proof. $\hfill \blacksquare$
  
 \parag{Remark} As a consequence of the previous corollary we have the following characterization of the roots of the $k$-th Bernstein polynomial of a $[\alpha]$-primitive fresco with nilpotent order $k$:
 \begin{itemize}
 \item $-\alpha -m$ is a root of the $[\alpha]$-primitive fresco with nilpotent order $k$ if and only if there exists a $\tilde{\mathcal{A}}$-linear surjective map of $ \mathcal{F}$ to a rank $k$ theme $T_k$
 such its $k$-th Bernstein polynomial is $(x + \alpha +m)$.
 \end{itemize}

\section{A more precise result than \cite{[B.22]}}

 Here is the strengthened version of the main result in \cite{[B.22]} announced in the introduction, which uses the higher order Bernstein polynomials of the concerned  fresco. The reader may remark that the main result of \cite{[B.23]} is a precise converse of this theorem.

\begin{thm}\label{improve}
Let $\alpha \in ]0,1]$ and assume the hypothesis $H(\alpha, 1)$\footnote{Remind that this mean that $\exp(2i\pi\alpha)$ is not a root of the local monodromy of $f$ acting on the reduced cohomology of the Milnor's fiber of $f$ at any point outside the origin.} for  the germ at the origin in $\mathbb{C}^{n+1}$ of holomorphic function $\tilde{f} : (\mathbb{C}^{n+1},0) \to (\mathbb{C}, 0)$. Assume that
 $\omega$ in $\Omega^{n+1}_0$ has the following property:\\
  there exists an integer $h \in \mathbb{Z}$ and a form $\omega' \in \Omega^{n+1}_0$ such that $F^{\omega, \omega'}_h(\lambda) $ has  a pole of order
  $p \geq 1 $ at some  point $\xi \in -\alpha -\mathbb{N}$. \\
  Note $\xi_p = -\alpha - m$  the biggest among  such  numbers $\xi \in -\alpha - \mathbb{N}$  for any choice of $\omega'$ and $h \in \mathbb{Z}$. Then the $p$-th Bernstein polynomial of the fresco $\mathcal{F}_{f, \omega}$ has a root in $[-\alpha -m, -\alpha] \cap \mathbb{Z}$.
\end{thm}

\parag{Proof} Note $P := P_1P_2$ the annihilator of the class of $[\omega]$ in the  $[\alpha]$-primitive quotient 
$$\mathcal{F}^\alpha := \mathcal{F}_{f,\omega}\big/\big(\mathcal{F}_{f,\omega}\big)_{\not= \alpha} $$
 of the fresco   $\mathcal{F}_{f,\omega}:= \tilde{\mathcal{A}}[\omega]$ inside the (a,b)-module $H^{n+1}_0$ associated to $f$, where $P_2$ is the annihilator of $[\omega]$ in $\mathcal{F}^\alpha\big/S_{p-1}\big(\mathcal{F}^\alpha\big)$. If $F^{\omega,\omega'}_h(\lambda)$ has a pole of order at least equal to $p$ at the point $-\alpha -m$ and if $-\alpha-m$ is not a root of the $q$-th Bernstein polynomial of $\mathcal{F}^\alpha$ for each $q \geq p$, then $-\alpha -m$ is not a root of the (usual) Bernstein polynomial of the fresco $\tilde{\mathcal{A}}\big/\tilde{\mathcal{A}}P_1$ which is isomorphic to $S_{p-1}(\mathcal{F}^\alpha)$. In this situation, using Corollary 8.2.5 in \cite{[B.23]} we see  that  $F^{P_2\omega, \omega'}_{h+p_2}(\lambda)$ has a pole of order at least equal to $p$ at $-\alpha -m$ , where $k$ is the rank of the fresco $\mathcal{F}^\alpha\big/S_{p-1}\big(\mathcal{F}^\alpha\big)$. But this is impossible, according to Corollary 8.2.6 in {\it loc. cit.} since the nilpotent order of  $S_{p-1}\big(\mathcal{F}^\alpha\big)$ is $p-1$ and since the image of the class of  $P_2\omega$ in $\mathcal{F}^\alpha$ generates $S_{p-1}(\mathcal{F}^\alpha)$.\\
   So there exists  integers $j \geq 0$ and $m_j \in [0,m]$ such that the  $(p+j)$-th Bernstein polynomial of $\mathcal{F}^\alpha_{f,\omega}$ has the root $-\alpha - m_j$. Then the $p$-th Bernstein polynomial of $\mathcal{F}_{f, \omega}$ has a root in $[-\alpha -m_j, -\alpha ] \subset [-\alpha - m, -\alpha]$ (thanks to the remark following Lemma \ref{missing}). $\hfill \blacksquare$\\

The end of Theorem 3.1.2 in \cite{[B.22]} is also improved as follows:

\begin{cor}\label{compl. 2} In the situation of the previous theorem, let, for  each integer $s$  in $[1,p]$,  $\xi_s$ be the biggest element in $-\alpha -\mathbb{N}$ for which there exists $h \in \mathbb{Z}$ and $\omega' \in \Omega^{n+1}_0$ such that $F^{\omega, \omega'}_h(\lambda)$ has a pole of order  at least equal to $s$ at $\xi_s$. Then $\xi_s$ is a root of  some  $(s+j)$-th Bernstein polynomial  of the fresco $\mathcal{F}_{f, \omega}^\alpha$ for some $j \in \mathbb{N}$.\\
Moreover, if $\xi_s = \xi_{s+1} = \cdots = \xi_{s+p}$, then there exists at least $p$ distinct values of $j \in \mathbb{N}$ such that $\xi_s$ is root of the $(s+j)$-th Bernstein polynomial of the fresco $\mathcal{F}_{f, \omega}^\alpha$.
\end{cor}

\parag{Proof} the proof of the first assertion is analogous to the proof of the theorem above. \\
The second assertion is an immediate consequence of the fact that the roots of the Bernstein polynomial of a semi-simple fresco are simple, applied to the successive semi-simple quotients
$$ S_d(\mathcal{F}^\alpha)\big/S_{d-1}(\mathcal{F}^\alpha) $$
for $d =  s+1,  s+2, \dots, s+p$.$\hfill \blacksquare$

\parag{Remarks}
\begin{enumerate}
\item If for $\omega$ given in $\Omega^{n+1}_0$, the maximal order of a pole at some point in $-\alpha - \mathbb{N}$ is equal to $d$ for any choice of $\omega' \in \Omega^{n+1}_0$ and any $h \in \mathbb{Z}$, then $\xi_d$ is a root of the $d$-th Bernstein polynomial of the fresco $\mathcal{F}_{f, \omega}$, because of the ``converse theorem"  proved in \cite{[B.23]} (see Theorem 8.5.1).

\item To consider a form  $\psi \in \mathscr{C}^\infty_c(\mathbb{C}^{n+1})^{0, n+1}$  with small enough support and such that $d\psi = 0$ in a neighborhood of $0$ is equivalent to consider $\rho\bar \omega'$ where $\omega' $ is in $\Omega^{n+1}_0$ and $\rho$ is a function in $\mathscr{C}^\infty_c(\mathbb{C}^{n+1})$ with small enough support which is identically $1$ near the origin. \\
Indeed any such $\psi$ may be written as $\psi = \bar \omega'$  for some $\omega' \in \Omega^{n+1}$ near the origin  thanks to Dolbeault's Lemma, and then $\psi -\rho \bar \omega'$ is identically $0$ near the origin, so replacing $\psi$ by $\rho \bar \omega' $ do not change the poles which may appear in $-\alpha - \mathbb{N}$  for the functions we are looking at (what ever is the choice of $h \in \mathbb{Z}$ thanks to our hypothesis $H(\alpha, 1)$).
\end{enumerate}

\section{ Examples}

 It is, in general, rather difficult to compute the Bernstein of the fresco associated to a given pair $(f, \omega)$, even in the case where $f$ has an isolated singularity. Nevertheless, in the case where $f$ is a polynomial in $\mathbb{C}[x_0, \dots, x_n]$ having $(n+2)$ monomials, we describe  in the article  \cite{[B.22]},  a rather elementary method to obtain an estimation for the Bernstein polynomial of the fresco $\mathcal{F}_{f,\, \omega}$ associated to a monomial $(n+1)$-form $\omega$.\\
Of course, when the full Bernstein polynomial has a  root of multiplicity $k \geq 2$ then this root is also a root of the $j$-th Bernstein polynomial for each $j \in [1,k]$ but when the Bernstein polynomial has only simple roots, the computation of the higher order Bernstein polynomials, even in the special situation of   \cite{[B.22]}, is not easy. We present in below some examples where we show that the second Bernstein polynomial  is not trivial but where the full  Bernstein polynomial has no multiple root.

\begin{prop}\label{example 0}
Let $f(x,y,z) := xy^3 + yz^3 + zx^3 + \lambda xyz $ where $\lambda \not= 0$ is any complex number  which is a parameter, and consider the holomorphic forms
 $$\omega_1 := dx\wedge dy\wedge dz,\  \omega_2 := y^3z^2 \omega_1, \ \omega_3 =  y^7 \omega_1, \ {\rm and}  \   \omega_4 := xy^3 \omega_1.$$
 Then, in each of these cases, the fresco $\mathcal{F}_{f, \,\omega_i}$ is a rank $2$ theme and the second Bernstein polynomial  is equal  respectively to $x+1$, $x+4$, $x+5$ and $x+3$.\\
  Moreover, for $i = 3, 4$ the corresponding (full) Bernstein polynomial of the corresponding frescos  has only simple roots.
 \end{prop}

Note that this proposition allows to apply Theorem 8.5.1 in \cite{[B.23]}  to conclude that for each $i \in \{1,2,3,4\}$, there exists some integer $h$ and some germ $\omega' \in \Omega^3_0$ such that  the meromorphic extension of 
$$ F_h^{\omega_i,\omega'}(\lambda) = \frac{1}{\Gamma(\lambda)}\int_X \vert f\vert^{2\lambda} \bar f^{-h} \rho \omega_i\wedge \bar\omega' $$
has a double pole at  the point $ \lambda_i$ equal to the root of the second Bernstein polynomial of the fresco $\mathcal{F}_{f,\omega_i}$.\\

The proof of this proposition uses several lemmas and the technic of computation described in \cite{[B.22]} (see paragraph 4.3.2) .

\begin{lemma}\label{gen. theme}
Let $e$ be a generator of the rank $2$ theme $T := \tilde{\mathcal{A}}/\tilde{\mathcal{A}}(a-2b)(a-b)$ (which is the unique fresco  with Bernstein polynomial $(x+1)^2$). Assume that we have three  homogeneous polynomials 
$P, Q$ and $R$ in $\mathcal{A}$ of  respective degrees $3,4$ and $k$ with the following conditions
\begin{enumerate}
\item  $P,Q $ and $R$ are monic in $a$.
\item Then exists a non zero constant $c$ such that $P + cQ$ kills $e$ in $T$.
\item The Bernstein polynomial of $Q$\footnote{By definition $B_P$ is defined by the formula $$ (-b)^pB_P(-b^{-1}a) = P$$ where $P$ is in $\mathcal{A}$, is homogeneous in (a,b) of degree $p$ and monic in $a$. This is the Bernstein polynomial of the fresco  $\tilde{\mathcal{A}}/\tilde{\mathcal{A}}P$.} is not a multiple of $(x+1)$ or of  $(x+2)$.
\item The Bernstein polynomial of $R$ is not a multiple of $(x+3)(x+2)(x+1)$
\end{enumerate}
Then $Re$ generates a rank two sub-theme in $T$.
\end{lemma}

\parag{Proof} First, remark that our hypothesis implies that $P = (a - \nu b)(a - 2b)(a - b)$ for some $\nu \in \mathbb{C}$ since $T$ is isomorphic to $\tilde{\mathcal{A}}\big/\tilde{\mathcal{A}}(a-2b)(a-b) $. We may realize $T$ in the simple pole asymptotic expansion module with rank $2$ which is isomorphic to $T^\sharp$
$$\Theta :=\Xi_1^2\big/\Xi_1^0 \simeq \mathbb{C}[[s]](Log\, s)^2 \oplus \mathbb{C}[[s]Log\, s$$
 where $a$ is the multiplication by $s$ and $b$ is defined by  $ab - ba = b^2$ and 
$$ b(Log\, s)= sLog\, s \quad {\rm and} \quad b((Log\,s)^2) = s(Log\, s)^2 - 2sLog\, s .$$ 
Then let us prove that  image of $e$ in $\Xi_1^2\big/\Xi_1^0 $ may be written 
\begin{equation*}
 e = u(Log\, s)^2 + vs(Log\, s)^2 + ws^3(Log\, s)^2 +  s^4\mathbb{C}[[s]] (Log\, s)^2+ \mathbb{C}[[s]](Log\, s) \tag{@}
\end{equation*}
where $\varphi$ is in $\Xi_1^2\big/\Xi_1^0 $ and where $uvw\not= 0$ are complex numbers. \\
  Remark that  the only restrictive condition for writing $e$ as in $(@)$ is the condition  $uvw\not= 0$. The condition $u \not= 0$ is easy because we assume that $e$ is a generator of $T$ with Bernstein polynomial $(x+1)^2$, so writing $e$ as a $\mathbb{C}[[b]]$-linear combination of the $\mathbb{C}[[b]]$-basis $e_1 = (Log\, s)^2$ and $e_2 = Log\, s$ of $T$  we see that the coefficient of $e_1$ must be invertible in $\mathbb{C}[[b]]$.\\
 But the condition $(P + cQ)(e) = 0$ implies, since the Bernstein element of $T$ is $(a-2b)(a-b)$, that we  may write\footnote{In our choice of $f$ and $\omega_1$, $\mu = 3$.} $P = (a-\nu b)(a - 2b)(a - b)$. \\
 The annihilator of  $(Log\, s)^2$ in $\Xi_1^2\big/\Xi_1^0 $ is the ideal $ \tilde{\mathcal{A}}(a-2b)(a-b)$ so we have $P( (Log\, s)^2 ) = 0$ in $T$. Since  $Q((Log\, s)^2)$ has a non zero term in $s^4(Log\, s)^2$, because $-1$ is not a root of $B_Q$, only the term coming from
 $$ P(s(Log\, s)^2) = \frac{4-\nu}{24}s^4(Log\, s)^2 \quad {\rm modulo} \ \mathbb{C}[[s]]Log\, s $$
 can compensate for this term, in order  to obtain  the equality $(P + cQ)(e) = 0$. Then $u\not= 0$. implies  $v \not= 0$. \\
 But now, the only term which can kill the non zero term in $s^5(Log\, s)^2$ coming from $Q(vs(Log\, s)^2)$ (using  that $B_{Q}$ is not a multiple of $(x+2)$)  can only comes from $P(ws^2(Log\,s)^2)$ and this proves that $w \not= 0$. So the assertion  $(@)$ holds true.\\
 Now if $R$ is homogeneous of degree $k$ in $(a,b)$ a necessary condition on $R$ such that $Re$ has no term in $s^{k+i}(Log\, s)^2$,  for $i= 0,1,2$,  is that $B_R$ divides $(x+1)(x+2)(x+3)$. So, when it is not the case Lemma  5.2.4 in \cite{[B.23]}  implies that $Re$ is a rank $2$ theme and its second Bernstein polynomial has a (unique)  root equal to $-(k+j)$ where $-j$ is the smallest integer among  $\{-1,-2,-3\}$ which is not a root of $B_R$ (see Corollary \ref{new 2}). $\hfill \blacksquare$\\
 
 Note that the Lemma above may be easily generalized  to many $[\alpha]$-primitive frescos  provided that the nilpotent order is known and that it has  a  generator which admits a enough simple element in its annihilator.\\
   
  \begin{lemma}\label{example 1}
 In the situation of Proposition \ref{example 0}, the frescos generated by the forms 
 $$\omega_1 := dx\wedge dy\wedge dz,  \ \omega_2 := y^3z^2 \omega_1, \ \omega_3 := y^7\omega_1, \  {\rm and} \  \omega_4 := xy^3\omega_1 $$ 
 generate rank $2$ $[1]$-primitive themes. Their Bernstein polynomials are respectively equal to $$(x+1)^2, \quad (x+3)^2 \ {\rm or} \ (x+2)(x+3), \quad  (x+3)(x+5) \quad {\rm and} \quad (x+2)(x+3)$$
 and  their respective $2$-Bernstein polynomials are  $(x+1)$, $(x+3)$, $(x+5)$ and $(x+3)$. In the  cases $i = 3,4$  there is no double root for the Bernstein polynomial of $\mathcal{F}_{f, \omega_i}.$ 
 \end{lemma}

 \parag{Proof} The first point is to show that $\mathcal{F}_{f,\omega_1}$ has  rank $2$. Since $f$ has an isolated singularity at the origin, we have $Ker df^{n}= df\wedge \Omega^{n-1}$ and then 
 $H^{n+1}/bH^{n+1} \simeq \mathcal{O}_{0}/J(f)$ and $H^{n+1}$ has no $b$-torsion and no $a$-torsion. Since $f$ is not\footnote{This point is not so easy to check directly. But  the rank is not $1$ since this would implies that  this fresco  has a simple pole and the argument used in  Lemma  \ref{gen. theme} gives then a contradiction.} 
   in $J(f)$ the image of $\omega_1$ and $a\omega_1 = f\omega_1$ in $H^{n+1}$ are linearly independent (over $\mathbb{C}$) and then the rank of $\tilde{\mathcal{A}}\omega_1$ is at least equal to $2$. Now the computation in \cite{[B.22]} (see 4.3.2) shows that the Bernstein polynomial of this fresco divides $(x+1)^3$ (see also the detailed computation below). So it is a theme of rank $2$ or $3$. But using our main result, the rank $3$ would imply that there exists a pole of order $3$ for some $F^{\omega_1,\omega'}_h(\lambda)$ which is impossible\footnote{This would give an order $4$ pole for the meromorphic continuation of $\vert f\vert ^{2\lambda}$ !} in $\mathbb{C}^3$. So $\mathcal{F}_{f,\omega_1}$ is a rank $2$ theme with Bernstein polynomial $(x+1)^2$. The computation in \cite{[B.22]} gives that $P_3 +  c\lambda^{-4}P_4$ kills $\omega_1$ in $H^{n+1}$ where 
 $$ P_3 := (a-3b)(a-2b)(a-b), P_4 = (a -(13/4)b)(a-(5/2)b)(a-(7/4)b)a, {\rm and } \ c = 4^4 $$
 This is easily obtain by using the technic of the computation of {\it loc.cit.} (see  the  detailed computation in the Appendix below). 
 Then we may apply  Lemma \ref{gen. theme} to see that $\lambda m_1m_2 \omega_1 = \lambda (a-2b)(a-b)\omega_1$  generates rank $2$ themes in $H^{n+1}$.
 But the identity $\lambda m_1m_2 = m_4 y^3z^2$ shows that $\omega_2$ generates also  rank $2$ in $H^{n+1}$ since $m_4\omega_2 = \lambda m_1m_2\omega_1 = \lambda(a-2b)(a-b)\omega_1$ applying Lemma  \ref{gen. theme} with $R = (a-2b)(a-b)$ whose Bernstein polynomial is $(x+1)^2$. Moreover we see that $Re$ has a non zero term in $s^3(Log\, s)^2$.\\
 Since $m_4\omega_2$ generates a rank $2$ theme, then $\omega_2$ generates a rank $2$ theme also (the rank $3$ is again excluded because  it would imply that $f^2 \not\in J(f)$ which is impossible as explained above).\\
 The technic of computation in \cite{[B.22]}  applied to $\omega_2$ gives now that the Bernstein polynomial of the rank $2$ theme $\tilde{\mathcal{A}}\omega_2$ has to divide\footnote{This computation gives that $Q_3 + d\lambda^{-4}Q_4$ kills $\omega_2$ in $H^{n+1}$ with $Q_3 := (a-4b)(a-4b)(a-3b)$.} the polynomial $(x+2)(x+3)^2$.\\
 But the fact that $m_4\omega_2$ has a  non zero term in $s^3(Log\, s)^2$ (and no term in $(Log\, s)^2$ or in $s(Log\,s)^2$)  implies, since we have 
  $$m_4\omega_2 = 4(a - 2b)\omega_2$$
  $\omega_2$ has a non zero term in  $s^2(Log\, s)^2$  and then $-3$ is a root of the second Bernstein polynomial of the fresco $\mathcal{F}_{f, \omega_2}$. So the Bernstein polynomial is either $(x+2)(x+3)$ or $(x + 3)^2$.\\
   We know\footnote{see the computation below.} that   the Bernstein polynomial of $\mathcal{F}_{f,\omega_3}$ divides $(x+5)(x+3)(x+2)$. 
 But we know also  that $m_1^2m_4\omega_1 = \lambda m_3\omega_3$ has a non zero term in  $s^5(Log\, s)^2$ (as a consequence of Lemma \ref{gen. theme}) and $-m_3\omega_3 = (a - 2b)\omega_3$ implies that $\omega_3$ has a non zero term in $s^4(Log\, s)^2$. Then  the second Bernstein polynomial of $\mathcal{F}_{f, \omega_3}$ is $x+5$.\\
  Note that, in this case, the Bernstein polynomial of the fresco $\mathcal{F}_{f, \omega_3}$ has two simple roots.\\
  The last  case is similar, since we know that $m_1\omega_1$  has a non zero term in $s^2(Log\, s)^2$.  So our assertion is consequence of the estimation of the Bernstein polynomial.
   $\hfill \blacksquare$\\
   
   For the convenience of the reader, we give in the Appndix below some detailed computations for these four examples.

\newpage

\section{ Appendix: detailed computations}
 
 \parag{The detailed computation for $\omega_1$}
  The linear system to compute $m_1, \dots, m_4$ (we note $\omega_1 $ by $1$) is given by $MX = \, ^t(a,b,b,b)(1)$ which gives
 \begin{align*}
 & 0 = (2m_1 +m_2 -3m_3) \quad {\rm and}\\
 & 0 = (3m_1 -2m_2 -m_3) \quad {\rm and \ so}\\
 & m_1 = m_2 = m_3  \quad {\rm and} \\
 &  m_4 = b(1) - 3m_1 -m_2= b(1) - 4m_1\\
 & a(1) = 3m_1 + b(1) -4m_1 \quad {\rm so} \quad m_1 = -(a-b)(1) = m_2 = m_3 \quad {\rm and} \quad m_4 = -(4a-3b)(1)\\
 \end{align*}
  The linear system to compute $m_1^2, \dots, m_1m_4$ is given by $MX = \, ^t(a,2b,4b,b)(m_1)$ which gives
 \begin{align*}
 & 2b(m_1) = (2m_1 +m_2 -3m_3)m_1 \quad {\rm and}\\
 & 3b(m_1) = (3m_1 -2m_2 -m_3)m_1 \quad {\rm and \ so}\\
 & m_1^2 = m_1m_2+ b(m_1) \quad m_1m_3 = m_1m_2  \quad {\rm and} \quad   m_1m_4 = b(m_1) - 4m_1m_2 \\
 & a(m_1) = m_1m_2 + b(m_1) + m_1m_2 + m_1m_2 +   b(m_1) - 4m_1m_2 \quad {\rm so} \\
 &  -m_1m_2 = (a -2b)(m_1) \quad {\rm and} \quad  -m_1^2 = (a- 3b)(m_1)= -(a-3b)(a-b)(1) \\
 \end{align*}
 The linear system to compute $m_1^2m_2,\  m_1m_2^2, \ m_1m_2m_3, \  m_1m_2m_4$ is given by  \\ $MX = \, ^t(a,2b,5b,4b)(m_1m_2)$ which gives
  \begin{align*}
 & 3b(m_1m_2) = (2m_1^2m_2 + m_1m_2^2 -3m_1m_2m_3) \quad {\rm and}\\
 & b(m_1m_2) = (3m_1^2m_2 -2m_1m_2^2 - m_1m_2m_3) \quad {\rm and \ so}\\
 &  m_1^2m_2 = m_1m_2^2, m_1m_2m_3 = m_1^2m_2 - b(m_1m_2), m_1m_2m_4 = 5b(m_1m_2) - 4m_1^2m_2 \\
 & -m_1^2m_2 = (a - 4b)(m_1m_2) \quad {\rm and} \quad -m_1m_2m_3 = (a- 3b)(m_1m_2). \\
 \end{align*}
 So the Bernstein polynomial of the fresco $\mathcal{F}_{f, \omega_1}$ has to divide the Bernstein polynomial of $P := (a-3b)(a-2b)(a-b)$ which is $B_P = (x+1)^3$. \\
 
  The linear system to compute $m_1m_4^p, \dots, m_4^{p+1}$ is given by  \\ $MX = \, ^t(a,(p+1)b,(p+1)b,(p+1)b)(m_4^p)$ which gives:
   \begin{align*}
 & 0 = (2m_1 +m_2 -3m_3)m_4^p \quad {\rm and}\\
 & 0  = (3m_1 -2m_2 -m_3)m_4^p  \quad {\rm and \ so}\\
 & m_1m_4^p = m_2m_4^p = m_3m_4^p  \quad {\rm and} \\
 &  m_4^{p+1} = (p+1)b(m_4^p) - (3m_1 +m_2)m_4^p = (p+1)b(m_4^p) - 4m_1m_4^p \quad {\rm so} \\
 & m_4^{p+1} = (4a - 3(p+1)b)(m_4^p) = 4(a - \frac{3}{4}(p+1)b)(m_4^p).\\
 \end{align*}
 So we obtain that $$Q = 4^4(a - 3b)(a - \frac{9}{4}b)(a - \frac{3}{2}b)(a - \frac{3}{4}b) = 4^4(a - \frac{13}{4}b)(a - \frac{5}{2}b)(a - \frac{7}{4}b)a.$$
 So $B_Q = x(x+1/4)(x+1/2)(x+3/4)$.
 
   \parag{The detail computation for $\omega_2 := y^3z^2\omega_1$} 
 Define  $\mu := y^3z^2$. The linear system to compute $m_1\mu, \dots, m_4\mu$ is given by $MX = \, ^t(a,b,4b,3b)(\mu)$ which gives
 \begin{align*}
 & 3b(\mu) = (2m_1 +m_2 -3m_3)\mu \quad {\rm and}\\
 & b(\mu) = (3m_1 -2m_2 -m_3)\mu \quad {\rm and \ so}
 & m_1\mu = m_2\mu = b(\mu) + m_3\mu  \quad {\rm and} \\
 &  m_4\mu = 4b(\mu) - 4m_1\mu 
 \end{align*}
 Then we obtain
 $$ -m_1\mu = -m_2\mu = (a -3b)(\mu), \quad -m_3\mu = (a-2b)(\mu)$$
 and also $ m_4\mu = 4(a-2b)(\mu)$  which is also  used above.\\
 Then we compute $m_1^2\mu, \dots, , m_1m_4\mu$ given by the system $MX = \, ^t(a,2b,7b,3b)(m_1\mu)$ which gives
 \begin{align*}
 & 5b(m_1\mu) = (2m_1 +m_2 -3m_3)m_1\mu \quad {\rm and}\\
 & 4b(m_1\mu) = (3m_1 -2m_2 -m_3)m_1\mu \quad {\rm and \ so}\\
 & (m_1 - m_2)m_1\mu = b(m_1\mu), \quad m_1m_3\mu = m_1m_2\mu - b(m_1\mu) \\
 & m_1m_4\mu = -4m_1m_2\mu + 4b(m_1\mu) 
 \end{align*}
 Then we obtain
 \begin{align*}
 &  -m_1m_2\mu = (a -4b)(m_1\mu), 
 \end{align*}
 Finally $m_1^2m_2\mu,  \ m_1m_2^2\mu, \  m_1m_2m_3\mu, \ m_1m_2m_4\mu$ are given by the system \\
   $MX = ^t(a,2b,8b,6b)(m_1m_2\mu)$ which gives
 \begin{align*}
 &  6b(m_1m_2\mu) = (2m_1 +m_2 -3m_3)m_1m_2\mu \quad {\rm and}\\
  & 2b(m_1m_2\mu) = (3m_1 -2m_2 -m_3)m_1m_2\mu \quad {\rm and \ so}\\
  & m_1^2m_2\mu = m_1m_2^2\mu, \quad  m_1m_2m_3\mu = m_1^2m_2\mu - 2b(m_1m_2\mu), \\
  & m_1m_2m_4\mu =  - 4m_1m_2m_3\mu.
  \end{align*}
  This gives $-m_1m_2m_3\mu = (a -4b)(m_1m_2\mu)$ and so the Bernstein element of the fresco $\mathcal{F}_{f, \omega_2}$ is given by
  \begin{equation*}
  m_1m_2m_3\mu = -(a-4b)(a-4b)(a-3b)(\mu) \tag{B}
  \end{equation*}
  which implies that the Bernstein polynomial divides $(x+2)(x+3)^2$. And we know that it is a rank $2$ theme with second Bernstein polynomial $x+3$. In this case it is not clear if the  Bernstein polynomial is $(x+3)^2$ or $(x+2)(x+3)$.
  $\hfill \square$
 
 \parag{The  detailed computation for $m_1^2m_4\omega_1 =  \lambda m_3\omega_3$}
 Since we have already obtained  $m_1^2 = (a - 3b)(a - b)(1)$ in the first case (computation for $\omega_1$),  it is enough to consider the linear system which  computes $m_1^3, \dots, m_4m_1^2$.\\
  It   is given by $MX = \, ^t(a,3b,7b,b)(m_1^2)$ and then
 \begin{align*}
 & 4b(m_1^2) = (2m_1 +m_2 -3m_3)m_1 \quad {\rm and}\\
 & 6b(m_1^2) = (3m_1 -2m_2 -m_3)m_1 \quad {\rm and \ so}\\
 & m_1^3 = m_1^2m_2+  2b(m_1^2) \quad m_1^2m_3 = m_1^2m_2  \quad {\rm and} \\
 &  m_1^2m_4 = b(m_1^2) - 4m_1^2m_2 \\
 & a(m_1^2) = m_1^2m_2 + 2b(m_1^2) + m_1^2m_2 + m_1^2m_2 +   b(m_1^2) - 4m_1^2m_2 \quad {\rm so} \quad  -m_1^2m_2 = (a -3b)(m_1^2) \\
 & m_1^2m_4 = b(m_1^2) +  4(a- 3b)(m_1^2) = (4a - 11b)(m_1^2)  
 \end{align*}
 So we obtain 
 $$ m_1^2m_4 = (4a - 11b)(a - 3b)(a - b)(1) .$$
 An the Bernstein polynomial of $R := (4a - 11b)(a - 3b)(a - b)$ is $(x+3/4)(x+2)(x+1)$ so we may apply Lemma \ref{gen. theme} and  $(@)$ gives that $R\omega_1 = \lambda m_3\omega_3$ has a non zero term in $s^5(Log\, s)^2$.

\parag{The  detailed computation for $ \omega_3 = y^7\omega_1$}
The linear system to compute $m_1\nu, \dots, m_4\nu$,  (where $\nu := y^7$) is given by $MX = \, ^t(a,b,8b,b)(y^7)$ which gives:
 \begin{align*}
 & 7b(\nu) = (2m_1 +m_2 -3m_3)\nu \quad {\rm and}\\
 & 7b(\nu) = (3m_1 -2m_2 -m_3)\nu \quad {\rm and \ so}\\
 & m_1\nu = m_2\nu + 2b(\nu), \quad  m_3\nu =  m_2 - b(\nu)  \quad {\rm and} \\
 &  m_4\nu  =  -4m_2\nu + 2b(\nu), \quad {\rm so} \quad -m_1\nu = (a - 5b)(\nu), \quad -m_2\nu = (a - 3b)(\nu), \\
 &  -m_3\nu = (a - 2b)(\nu) \quad {\rm and} \quad m_4 = (4a-10b)\\
 \end{align*}
 
 Remark that this computation is already enough to see that the Bernstein polynomial of $\mathcal{F}_{f,\,\omega_3}$ divides $(x+5)(x+3)(x+2)$\footnote{because we can compute $m_1m_2m_3\nu$ with any order for $1, 2, 3$ and then obtain that the $P$ which estimates the Bernstein element of the fresco $\mathcal{F}_{f, \omega_3}$ is right divisible by $(a - 5b), (a-3)$ and $(a-2b)$.}. 
 But we know that $m_1^2m_4\omega_1 = \lambda m_3\omega_3$ has a non zero term in  $s^5(Log\, s)^2$ and $-m_3\omega_3 = (a - 2b)\omega_3$ implies that $\omega_3$ has a non zero term in $s^4(Log\, s)^2$ and then  the second Bernstein polynomial of $\mathcal{F}_{f,\omega_3}$ is $x+5$. The Bernstein polynomial is $(x+3)(x+5)$ because $(a-2b)\omega_3$ cannot have a non zero term in $s(Log\, s)^2$ (note that  $m_1^2m_4\omega_1$ has a non zero term in $s^3Log\, s$)  since $u\not= 0$ and so $\omega_3$ has a non zero term in $s^2Log\, s$).

   \parag{The  detailed computation for $\omega_4 := m_1\omega_1$}
 The linear system to compute $m_1^2, \dots, m_1m_4^2$ is given by $MX = \, ^t(a,2b,4b,b)(m_1)$ is already solved above and it gives
 $$  m_1^2 = (a - 3b)(m_1).$$

 The linear system to compute $m_1^2m_2, \dots, m_1^2m_4$ is given by $MX = \, ^t(a,3b,7b,b)(m_1^2)$ which gives:
  \begin{align*}
  & 4b(m_1^2) = (2m_1 +m_2 -3m_3)m_1m_4 \quad {\rm and}\\
 & 6b(m_1^2)  = (3m_1 -2m_2 -m_3)m_1m_4  \quad {\rm and \ so}\\
 & m_1^3 = m_1^2m_2 + 2b(m_1^2), \quad {\rm and} \  m_1^2m_3 = m_1^2m_2   \quad {\rm and} \\
 &  m_1^2m_4 = 7b(m_1m_4) - 3m_1^3 -m_1^2m_2 = b(m_1^2) - 4m_1^2m_2 \quad {\rm so} \quad -m_1^2m_2 = (a - 3b)(m_1^2).
 \end{align*}
  The linear system to compute $m_1^3m_2, m_1^2m_2^2, m_1^2m_2m_3, m_1^2m_2m_4$ is given by  \\
  $MX = \, ^t(a,3b,8b,4b)(m_1^2m_2)$ which gives:
  \begin{align*}
  & 5b(m_1^2m_2) = (2m_1 +m_2 -3m_3)m_1^2m_2 \quad {\rm and}\\
 & 4b(m_1)  = (3m_1 -2m_2 -m_3)m_1^2m_2  \quad {\rm and \ so}\\
 & m_1^3m_2 = m_1^2m_2^2 + b(m_1^2m_2), \ m_1^2m_2m_3 =  m_1^2m_2^2 - b(m_1^2m_2)  \quad {\rm and} \\
 &  m_1^2m_2m_4  = 8b(m_1) - 3m_1^2m_2 - m_1^2m_2^2= 5b(m_1^2m_2) - 4m_1^2m_2^2 \quad {\rm so} \\
 &  m_1^2m_2m_3 = (a - 4b)(m_1^2m_2).
 \end{align*}
 Finally we find $(a - 4b)(a - 3b)(a - 3b)(m_1\omega_1) = -m_1m_2m_3(m_1\omega_1)$. So the Bernstein polynomial of the fresco $\mathcal{F}_{f, \, \omega_4}$ has to divide $(x+2)^2(x+3)$. Since $m_1\omega_1 = -(a - b)\omega_1$ has clearly a non zero term in $s^2(Log\, s)^2$ this fresco is rank $2$ theme and its second Bernstein polynomial is $x+3$. Then its Bernstein polynomial is $(x+2)(x+3)$.

\end{document}